\documentclass[12pt,amscd]{amsart}
\usepackage[all]{xy}
\usepackage{graphicx}
\usepackage{amsmath,amsxtra,amssymb,latexsym, amscd,amsthm}
\usepackage{indentfirst}
\usepackage[mathscr]{eucal}
\usepackage[T1]{fontenc}

\newtheorem{thm}{Theorem}[section]
\newtheorem{thm*}{Theorem}[]
\newtheorem{ques}{Question}[]
\newtheorem{cor}[thm]{Corollary}
\newtheorem{lem}[thm]{Lemma}
\newtheorem{prop}[thm]{Proposition}
\theoremstyle{definition}
\newtheorem{defn}[thm]{Definition}
\newtheorem{exm}[thm]{Example}

\newtheorem{conj}[thm*]{Conjecture}

\DeclareMathOperator{\N}{\mathbb {N}}
\DeclareMathOperator{\Z}{\mathbb {Z}}
\DeclareMathOperator{\R}{\mathbb {R}}
\DeclareMathOperator{\depth}{depth}
\DeclareMathOperator{\ass}{ass}
\DeclareMathOperator{\spec}{Spec}

\DeclareMathOperator{\supp}{supp}

\DeclareMathOperator{\asr}{assrad}
\DeclareMathOperator{\Min}{Min}
\DeclareMathOperator{\bht}{bight}

\newcommand{\pp}{\mathfrak p}
\def\alb {\boldsymbol {\alpha}}
\def\btb {\boldsymbol {\beta}}

\def\gmb {\boldsymbol {\gamma}}
\def\lam {\boldsymbol {\lambda}}

\def\e {\mathbf e}
\def\x {\mathbf x}
\def\x {\mathbf x}

\def\u {\mathbf u}

\def\p {\mathfrak p}

\def\H {\mathcal{H}}
\def\V {\mathcal{V}}
\def\E {\mathcal{E}}

\begin{document}

\title[On the set of associated radicals of powers of monomial ideals] {On the set of associated radicals of powers of monomial ideals}

\author[N.T. Hang]{Nguyen Thu Hang}
\address{Thai Nguyen University of Sciences, Thai Nguyen, Vietnam}
\email{hangnt@tnus.edu.vn}

\author[T.T. Hien]{ Truong Thi Hien}
\address{Hong Duc University, Thanh Hoa, Viet Nam, and Institue of Mathematics, VAST, Hanoi, Viet Nam}
\email{hientruong86@gmail.com}

\date{}

\dedicatory{}
\commby{}
%-----------------------------------------------------------
\begin{abstract} Let $I$ be a monomial ideal in a polynomial ring. In this paper, we study the asymptotic behavior of the set of associated radical ideals of the (symbolic) powers of $I$. We show that both $\asr(I^s)$ and $\asr(I^{(s)})$ need not stabilize for large value of $s$. In the case $I$ is a square-free monomial ideal, we prove that $\asr(I^{(s)})$ is constant for $s$ large enough. Finally, if $I$ is the cover ideal of a balanced hypergraph, then $\asr(I^s)$ monotonically increases in $s$. 
\end{abstract}
% -----------------------------------------------------------
\maketitle
% -----------------------------------------------------------
\section{Introduction}
Let $R = K[x_1,\ldots,x_n]$ be a polynomial ring over a field $K$ and  $I$ a monomial ideal of $R$. The study of $\depth R/I$ in terms of some combinatorial structure associated with $I$ is an interesting problem in combinatorial commutative algebra (see \cite{FHM, MV}). One of the efficient ways to compute this invariant is  by using Hochster's formula (see \cite{HOC}):
$$\depth R/I = \min\{\depth R/\sqrt{I\colon f} \mid f \text{ is a monomial not belong to } I\}.$$
This formula is particularly useful when $I$ is not a square-free monomial ideal (see e.g. \cite{BJ, HHV, MTT1, MTT2}).

A  square-free monomial ideal $J$ is called an associated radical ideal of $I$ if $J = \sqrt{I\colon f}$ for some monomial element $f\notin I$ (see \cite{HOC}). Let $\asr(I)$ denote the set of all associated radical ideals of $I$ (see \cite{BJ,JS}). Using Hochster's formula, we can deduce that the sequence $\{\asr(I^s)\}_{s\geqslant 1}$ is increasing (or stationary), and so is $\depth R/I^s$. A similar fact holds true for the symbolic powers of $I$. Thus, the combinatorial objects $\asr(I^s)$ and $\asr(I^{(s)})$ can be used to study $\depth R/I^s$ and $\depth R/I^{(s)}$, respectively.

This paper focuses on studying the asymptotic behavior of the set of associated radical ideals of (symbolic) powers of $I$. Let $\ass(I)$ be the set of associated prime ideals of $R/I$. Then, we have $\ass(I) \subseteq \asr(I)$ (see Lemma \ref{ass_assrad}).  Motivated by Brodmann's results \cite{B1, B2} which say that both $\ass(I^s)$ and $\depth R/I^s$ are constant for $s\gg 0$, it is natural to ask the following question for a monomial ideal $I$ of $R$.

\medskip

\begin{ques}\label{Quest1}
Is $\asr(I^s)$ constant in $s$ for $s\gg 0$? 
\end{ques}

\medskip

A similar question for the symbolic powers of $I$ is as follows.

\medskip

\begin{ques}\label{Quest2}
Is $\asr(I^{(s)})$ constant in $s$ for $s\gg 0$?    
\end{ques}

\medskip

However, the above two questions need not have positive answers in general. Namely,

\begin{thm*} [Theorem \ref{T1}] There is a monomial ideal $I$ in the polynomial ring $R = K[x,y,z,u,v]$ such that both $\asr I^{(s)}$ and $\asr (I^s)$ are not constant for $s\gg 0$.
\end{thm*}

Therefore, Questions \ref{Quest1} and \ref{Quest2} are now restricted to the special classes of monomial ideals. In the case $I$ is a square-free monomial ideal, the main result of the paper is the following theorem.

\begin{thm*} [Theorem \ref{T2}] Let $I$ be a square-free monomial ideal in $R$. Then, $\asr I^{(s)} = \asr I^{(s_0)}$ for all $s\geqslant s_0$, where $s_0=\lceil n\bht(I)^{(n+2)/2}\rceil$. In particular, $\asr I^{(s)}$ is constant for $s\gg 0$.
\end{thm*}

By using this theorem  we can prove that $\depth R/I^{(s)} = \dim R - \ell_s(I)$ for $s\geqslant n\bht(I)^{(n+2)/2}$ where $\ell_s(I)$ is the symbolic analytic spread of $I$ (see Corollary \ref{DETP-STABILITY}), it is a slight improvement of \cite[Theorem 2.4]{HKTT}.

\medskip

For the ordinary powers of $I$, we propose the following conjecture.

\begin{conj} Let $I$ be a square-free monomial ideal. Then, $\asr I^s$ is constant for $s\gg 0$.
\end{conj}

Finally, if the square-free monomial $I$ is the cover ideal of a balanced hypergraph, then the behavior of the sequence $\{\asr(I^s)\}_{s\geqslant 1}$ is nice. In fact, it is an increasing sequence. It is worth mentioning that this property is not preserved by any square-free monomial ideal (see Theorem \ref{not-increasing-sequence}).

\begin{thm*} [Theorem \ref{T3}] Let $I=J(\H)$ be a cover ideal of a balanced hypergraph $\H$. Then, $\asr (I^{s})\subseteq \asr (I^{s+1})$ for every $s\geqslant 1$. 
\end{thm*}

The paper is organized as follows. In Section $1$, we set up the notations, provide some background about the set of associated radical ideals of (symbolic) powers of a monomial ideal $I$, and show that $\asr(I^s) $ and $\asr(I^{(s)}) $ are not constant for $s\gg 0$. In Section $2$, we establish the assymptotic stability of $\asr(I^{(s)}) $ when $I$ is a square-free monomial ideal. In Section $3$, by investigating the vertices of the polytope, we prove the increasing property of $\asr(I^s) $ and $\asr(I^{(s)}) $ in the case $I$ is a cover ideal of a balanced hypergraph.

\section{Preliminaries}
In this section, for the reader's convenience, we recall some basic notations used in this paper and several auxiliary results.

Let $R=K[x_1,\ldots,x_n]$ be a polynomial ring over a field $K$. Let $I$ be a monomial ideal of $R$.  Then the $s$-th {\it symbolic power} of $I$, $s\geqslant 1,$ is defined by
$$I^{(s)}= \bigcap_{\pp\in \Min(I)} I^sR_\pp\cap R,$$
where $\Min(I)$ is  the set of minimal associated prime ideals of $I$.

We can compute the symbolic powers of monomial ideal $I$ by using its primary decompositions. Let
\begin{equation}\label{decomposition}
I = Q_1\cap \cdots\cap Q_r\cap Q_{r+1}\cap \cdots \cap Q_t
\end{equation}
be a minimal primary decomposition of $I$, where $Q_i$ is a primary monomial ideal for $i=1,\ldots,t$, and $\pp_j=\sqrt{Q_j}$ is a minimal prime of $I$ if and only if $1\le j \le r$ (hence $Q_{r+1},\ldots,Q_t$ are embedded primary components) so that $\Min(I) =\{\pp_1,\ldots,\pp_r\}$. Then, we have
\begin{equation}\label{decompositionprimary}
I^{(s)} = Q_1^s \cap Q_2^s\cap\cdots\cap Q_r^s.
\end{equation}

If $I$ is a square-free monomial ideal in $R$, then by \cite[Corollary 1.3.6]{HH1}, $I$ has a minimal primary decomposition of the form
\begin{equation}\label{square-free-dmp}
I = \p_1\cap \cdots\cap \p_r    
\end{equation}
and each $\pp_i\in \ass(I)$ is a monomial prime ideal that is generated by some variables. Hence, by \eqref{decompositionprimary}, for any $s\geqslant 1$, we have
\begin{equation}\label{decompositionprimary-dmp}
I^{(s)} =  \pp^s_1\cap \ldots \cap \pp^s_r.
\end{equation}

We also present a square-free monomial ideal via a hypergraph (see \cite{Berge} for more details on hypergraphs). Let $\V = \{1,\ldots,n\}$, and let $\E$ be a family of distinct non-empty subsets of $\V$. The pair $\H=(\V,\E)$ is called {\it a hypergraph} with the vertex set $\V$ and the edge set $\E$. Note that a hypergraph generalizes the classical notion of a graph. It means that a graph is a hypergraph for which every $E \in \E$ has cardinality two.

A {\it vertex cover} of $\H$ is a subset of $\V$ which intersects every edge of $\H$; a vertex cover is {\it minimal} if none of its proper subsets is itself a cover. For a subset $\tau = \{i_1,\ldots,i_t\}$ of $\V$, set $\x_{\tau} = x_{i_1}\cdots x_{i_t}$. The {\it cover ideal} of $\H$ is then defined by: 
$$J(\H) = (\x_{\tau} \mid \tau \text{ is a minimal vertex cover of } H).$$

It is well-known that there is a one-to-one correspondence between square-free monomial ideals of $R$ and cover ideals of hypergraphs on the vertex set $\V$.

The cover ideal of $J(\H)$ can be written in terms of primary decomposition as
\begin{equation} \label{intersect1}
J(\H) = \bigcap_{E\in\E} (x_i\mid i\in E).
\end{equation}
 
The following concept is the subject of this paper.

\begin{defn}\label{assradset} Let $I$ be a monomial ideal in $R$ and $f$ be a monomial, which is not in $I$. The radical ideal $\sqrt{I:f}$ is called an associated radical ideal of $I$. We denote the set of all associated radical ideals of I by $\asr (I)$.
\end{defn} 

With this definition, the Hochster's formula (see \cite[Theorem 7.1]{HOC}) is written as:
\begin{equation}\label{Hochsterfomular}
\depth R/I=\min \{\depth R/J \mid J\in \asr(I)\}.
\end{equation} 

\begin{exm} Let $I=(x_1,x_2^2)\cap (x_2,x_3^3)\cap (x_1,x_3)$, $P_1=(x_1,x_2), P_2=(x_2,x_3)$ and $P_3=(x_1,x_3)$. Then, $\asr (I)=\{P_1, P_2, P_3, P_1\cap P_2, P_2\cap P_3, P_1\cap P_2\cap P_3\}.$
\end{exm}

Denote $\ass(I)$ to be the set of associated prime ideals of $R/I$.

\begin{lem}\label{ass_assrad} $\ass(I) = \asr(I) \cap \spec(R)$.
\end{lem}
\begin{proof} If $\p\in\ass(I)$, then $\p = I \colon f$ for some monomial $f\notin I$. It follows that
$\p \in \asr(I) \cap \spec(R)$, so $\ass(I) \subseteq \asr(I) \cap \spec(R)$.

Conversely, for any $\p \in \asr(I) \cap \spec(R)$, assume $\p = \sqrt{I \colon f}$ for some monomial $f\notin I$. Let $I = Q_1\cap \cdots\cap Q_t$ be the minimal primary decomposition of $I$. Then, $\p=\sqrt{I\colon f} = \bigcap_{i=1}^t \sqrt{Q_i\colon f}$. It follows that $\p = \sqrt{Q_i\colon f}$ for some $i$. In particular, $\p=\sqrt{Q_i}\in\ass(I)$, thus $\asr(I) \cap \spec(R)\subseteq\ass(I)$, and hence the lemma follows.
\end{proof}

Let $G(I)$ be the set of minimal monomial generators of $I$. We define the support of $I$, denoted by $\supp(I)$, to be
$$\supp(I) = \{i\in [n] \mid  x_i \text{ divides } f \text{ for some } f \in G(I)\},$$
i.e., it is the set of the indices of variables appearing in generators of $I$.

In order to facilitate an induction argument on the number of variables we use the following notion (see \cite[Definition 8]{T}). For a vector $\alb = (\alpha_1,\ldots,\alpha_n)\in \N^n$ and for each $i= 1,\ldots, n$, let $\x^{\alb}[i]$ denote the monomial $x_1^{\alpha_1}\cdots \widehat{x^{\alpha_i}}\cdots x_n^{\alpha_n}$, where the term $x_i^{\alpha_i}$ of $\x^{\alb}$ is omitted. For a monomial ideal $I$ we define
$$I[i] = (\x^{\alb}[i]\mid \x^{\alb}\in I).$$
Actually, $I[i] = IR_{x_i} \cap R$, where $R_{x_i}$ is the localization of $R$ with respect to $x_i$.

Since $(\sqrt{I\colon J})[i] = \sqrt{I[i] \colon J[i]}$ for any other monomial ideal $J$ and for each $i$, we obtain the following fact.

\begin{lem}\label{localization} $\asr(I)=\bigcup_{i=1}^n \asr(I[i]) \cup \{J\in \asr(I)\mid \supp(J) = [n]\}$.
\end{lem}

Brodmann \cite{B1} showed that $\ass(I^s)$ is constant for $s\gg 0$. By contrast, both $\asr I^{(s)}$ and $\asr I^s$ may not be constant for $s\gg 0$. In order to prove this fact, we first prove the following lemma.

\begin{lem} \label{inc-lemma} If $\asr(I^r) \subseteq \asr(I^s)$ for some $r,s\geqslant 1$, then $\asr(I^{(r)}) \subseteq \asr(I^{(s)})$.
\end{lem}
\begin{proof} Let $J \in \asr(I^{(r)})$ so that $J = \sqrt{I^{(r)}\colon f}$ for some monomial $f\notin I^{(r)}$. Note that $I^r$ has a minimal primary decomposition of the form
$$I^r = I^{(r)} \cap Q_{1,r}\cap\cdots \cap Q_{u,r}$$
where $Q_{1,r},\ldots, Q_{u,r}$ are primary ideals of $I^r$ such that each $P_{i,r} = \sqrt{Q_{i,r}}\ , i= 1, \ldots, u$ is an embedded prime ideal of $I^r$.

Similarly, $I^s$ has a minimal primary decomposition of the form
$$I^s = I^{(s)} \cap Q_{1,s}\cap\cdots \cap Q_{v,s}$$
where $Q_{1,s},\ldots, Q_{v,s}$ are primary ideals of $I^s$ such that each $P_{j,s} = \sqrt{Q_{j,s}} \ , j = 1, \ldots, v$ is an embedded prime ideal of $I^s$.

Since $I^r \subset I^{(r)}$, we have $f\notin I^r$. Thus, $\sqrt{I^r\colon f}\in \asr(I^r)$. Since $$\sqrt{I^r:f}=\sqrt{(I^{(r)} \cap Q_{1,r}\cap \ldots \cap Q_{u,r}):f}=\sqrt{I^{(r)}:f}\cap \left(\bigcap_{i=1}^u \sqrt{Q_{i,r}\colon f}\right),$$
we have
\begin{equation}\label{FC1}
\sqrt{I^r:f}=J \cap \left(\bigcap_{f\notin Q_{i,r}} P_{i,r}\right).    
\end{equation}

By our assumption, $\asr(I^r)\subseteq\asr(I^s)$, so 
$$J \cap \left(\bigcap_{f\notin Q_{i,r}} P_i\right) = \sqrt{I^s\colon g}$$
for some $g\notin I^s$. By the same argument as in the proof of Formula $(\ref{FC1})$ we get
\begin{equation}\label{FC2}
\sqrt{I^s:g}=\sqrt{I^{(s)}\colon g} \cap \left(\bigcap_{g\notin Q_{j,s}} P_{j,s}\right).
\end{equation}

Note that both $J$ and $\sqrt{I^{(s)}:g}$ are intersections of some minimal prime ideals of $I$. From $(\ref{FC1})$ and $(\ref{FC2})$ we deduce that $\sqrt{I^{(s)}\colon g} = J$. Hence, $J\in\asr(I^{(s)})$, and the lemma follows.
\end{proof}

\begin{thm} \label{T1} There is a monomial ideal $I$ in the polynomial ring $R = K[x,y,z,u,v]$ such that both $\asr (I^{(s)})$ and $\asr (I^s)$ are not constant for $s\gg 0$.
\end{thm}
\begin{proof}  Let $I = (x^2,y^2,z^2)^2 \cap (x^3,y^3,u)\cap (z,u)$ be a monomial ideal of $R$. By \cite[Lemma 4.5]{NT}, we obtain
\begin{equation}\label{depthSymbolicPowers}
\depth R/I^{(s)}=
\begin{cases}
1 & \text{ if } s \text{ is odd},\\
2 & \text{ if } s \text{ is even}.
\end{cases}
\end{equation}

Together with Hochster's formula, it yields $\asr(I^{(s)})$ is not constant for $s\gg 0$. Together with Lemma \ref{inc-lemma}, it also deduces that $\asr (I^s)$ is not a constant for $s\gg 0$, and the theorem follows.
\end{proof}

\section{Symbolic powers of square-free monomial ideals}

In this section, let $I$ be a square-free monomial ideal in $R$. We will prove that $\asr (I^{(s)})$  is constant for $s\gg 0$. 

By Equation $(\ref{square-free-dmp})$ we may assume that $I$ has the primary decomposition given by
\begin{equation}\label{square-free-decomp}
I = \p_1\cap \cdots\cap \p_r, \text{ with } \ass(I) = \{\p_1,\ldots,p_r\},
\end{equation}
and for any $s\geqslant 1$, we have
\begin{equation}\label{symbolic-decompositionprimary}
I^{(s)} =  \pp^s_1\cap \ldots \cap \pp^s_r.
\end{equation}

Let $\x^{\alb}$ be a monomial of $R$ where $\alb=(\alpha_1,\ldots,\alpha_n)$. Then,
$$\sqrt{I^{(s)}\colon \x^{\alb}} = \bigcap_{\x^{\alb}\notin \pp_i^s} \pp_i.$$

Note that $\x^{\alb}\notin \p_i^s$ if and only if $$\sum_{x_j\in \pp_i} \alpha_j \leqslant s-1.$$

It yields the following useful fact.

\begin{lem}\label{FF} For any $\alb = (\alpha_1,\ldots,\alpha_n)\in \N^n$ we have
$$\sqrt{I^{(s)}\colon \x^{\alb}} = \bigcap \{\p \mid \p\in \ass(I) \text{ and } \sum_{x_j\in\p}\alpha_j \leqslant s-1\}.$$  
\end{lem}

\begin{lem} $\asr(I^{(s)}) \subseteq \asr(I^{(st)})$ for all $s,t\geqslant 1$.
\end{lem}
\begin{proof} Assume that $\ass(I)=\{\p_1,\ldots,\p_r\}$. Let $J= \sqrt{I^{(s)}\colon \x^{\alb}}\in \asr(I^{(s)})$. We also may assume that $J = \p_1\cap\cdots\cap \p_k$ for some $1\leqslant k \leqslant r$. By Lemma \ref{FF} we have the linear inequality
$$\begin{cases}
\sum\limits_{x_i\in \p_j} \alpha_i  \leqslant  s-1 & \text{ for } j = 1,\ldots,k,\\
\sum\limits_{x_i\in \p_j} \alpha_i  \geqslant  s & \text{ for } j = k+1,\ldots,r.
\end{cases}
$$
Consequently,
$$\begin{cases}
\sum\limits_{x_i\in \p_j} t\alpha_i  \leqslant  t(s-1)\leqslant ts-1 & \text{ for } j = 1,\ldots,k,\\
\sum\limits_{x_i\in \p_j} t\alpha_i  \geqslant  ts & \text{ for } j = k+1,\ldots,r.
\end{cases}
$$
Together with Lemma \ref{FF} again, it forces $\sqrt{I^{(st)}\colon \x^{t\alb}} = J$. Thus, $J\in\asr(I^{(st)})$, and the lemma follows.
\end{proof}

By virtue of this lemma, it is natural to ask whether $\asr(I^{(s)})$ is increasing in $s$, i.e., $\asr(I^{(s)}) \subseteq \asr(I^{(s+1)})$  for $s\geqslant 1$. Unfortunately, it is not true as shown in the following theorem.

\begin{thm}\label{not-increasing-sequence} There is a square-free monomial ideal $I$ in a some polynomial ring such that both $\asr(I^{(s)})$ and $\asr(I^s)$ are not increasing in $s$.
\end{thm}
\begin{proof} By \cite[Theorem 2.8]{NT}, for $s\geqslant 3$, there is a square-free monomial ideal in some polynomial ring $R$ such that $\depth R/I^{(s)} < \depth R/I^{(t)}$ for some $t > s$. Together with Hochster's formula, we deduce that $\asr(I^{(s)}) \not\subseteq \asr(I^{(t)})$, i.e.,  the sequence $\{\asr(I^{(s)})\}_{s\geqslant 1}$ is not an increasing sequence.

Together this fact with Lemma \ref{inc-lemma}, we deduce that $\asr(I^s)$ is also not increasing in $s$, as required.
\end{proof}

In the rest of this section, we will prove that the sequence $\{\asr(I^{(s)})\}_{s\geqslant 1}$ is stationary. Suppose that $I$ has the minimal primary decomposition given by (\ref{square-free-decomp}). For each associated radical ideal $Q$ in $\asr (I^{(s)})$, we have $Q=\sqrt{I^{(s)}:\x^{\alb}}$ for some monomial $\x^{\alb}=x_1^{\alpha_1}\ldots x_n^{\alpha_n}\notin I^{(s)}$. 

Let $F_i = \supp(\p_i)$ for $i=1,\ldots,r$. We may assume that $Q = \p_1\cap \cdots\cap \p_k$ for some $1\leqslant k\leqslant r$. By Lemma \ref{FF}, $\alb\in \Lambda_s$, where $\Lambda_s$ is a convex polyhedron in $\R^n$ defined by the  system of linear inequalities:
\begin{equation}\label{alphaRoot}
\begin{cases}
\sum\limits_{i\in F_j} x_i\leqslant s-1, \text{ for } j=1,\ldots,k,\\
\sum\limits_{i\in F_j} x_i\geqslant s, \text{ for } j=k+1,\ldots,r,\\
x_1\geqslant 0, \ldots,x_n\geqslant 0.
\end{cases}    
\end{equation}
Moreover, $Q=\sqrt{I^{(s)}\colon \x^{\btb}}$ for every $\btb\in \Lambda_s\cap \N^n$.

Hence, it is an interesting problem to investigate the integer solutions of this system. For this purpose, with each $t\geqslant 1$, we denote by $C_t $ to be the set of solutions  in  $\mathbb R^n$ of the following system of linear inequalities:
\begin{equation}\label{EQ-basics}
\begin{cases}
\sum\limits_{i\in F_j} x_i  <  t & \text{ for } j =1,\ldots,k,\\
\sum\limits_{i\in F_j} x_i  \geqslant  t & \text{ for } j =k+1,\ldots,r,\\
x_1\geqslant 0,\ldots,x_n\geqslant 0.
\end{cases}
\end{equation}

Then, $\alb \in C_s \cap \N^n$. Moreover, for $t\geqslant 1$, by Lemma \ref{FF} we have $Q\in\asr(I^{(t)})$ if and only if $C_t \cap \N^n \ne \emptyset$. Let $\overline{C_t}$ be the closure of $C_t$ in $\R^n$ with respect to the Euclidean topology. It means  $\overline{C_t}$ is  the set of solutions of the following system:
 \begin{equation}\label{EQ-C}
\begin{cases}
\sum\limits_{i\in F_j} x_i  \leqslant  t & \text{ for } j =1,\ldots,k,\\
\sum\limits_{i\in F_j} x_i  \geqslant  t & \text{ for } j =k+1,\ldots,r,\\
x_1\geqslant 0,\ldots,x_n\geqslant 0.
\end{cases}
\end{equation}

We now investigate the structure of $\overline{C_t}$. It turns out that it is a polytope in $\R^n$ if $\Lambda_s\ne\emptyset$.

\begin{lem}\label{polytope} Suppose that $\bigcup_{i=1}^k F_i = [n]$ and $\Lambda_s\ne \emptyset$ for some $s$. Then, $\overline{C_t}$ is a polytope in $\mathbb R^n$ for all $t\geq 1$, and  $\dim \overline{C_t}=n$.
\end{lem}

\begin{proof} Since $\overline{C_t} = t\overline{C_1}$, it suffices to prove the lemma for $\overline{C_1}$. First, we prove $C_1$ is bounded. Indeed, let $y=(y_1,\ldots,y_n)\in C_1$. For every $j\in [n]$, since $\bigcup_{i=1}^k F_i = [n]$, we have $j\in F_i$ for some $i\in \{1,\ldots,k\}$. Together with the system \eqref{EQ-C} where $t=1$, we get 
$$0\leqslant   y_j\leqslant \sum\limits_{v\in F_i} y_v<1,$$
so $C_1$ is bounded, as claimed.

It follows that $\overline{C_1}$ is bounded, so it is a polytope in $\R^n$. Now we prove $\dim \overline C_1 = n$. Note that $\Lambda_s\subseteq C_s$, so that $C_s\ne \emptyset$. It follows that $C_1\ne\emptyset$ since $C_t = tC_1$. Therefore, we can take a point $\btb = (\beta_1, \ldots , \beta_n) \in C_1$. From the system $(\ref{EQ-basics})$ with $t=1$, there exists a real number $\varepsilon > 0$ such that for all real numbers $\varepsilon_1, \ldots,\varepsilon_n$ with $0 \leqslant \varepsilon_1,\ldots,\varepsilon_n\leqslant \varepsilon$, and $\e_i=(0, \ldots,1, \ldots,0)$ where the number $1$ is in $i$-th position, we have $\btb+\varepsilon_1\e_1+\cdots+\varepsilon_n \e_n \in C_1$. Thus, $[\beta_1, \beta_1+\varepsilon] \times \cdots \times [\beta_n, \beta_n+\varepsilon]\subseteq C_1\subseteq \overline C_1$. This implies that the polytope $\overline C_1$ is full dimensional in $\R^n$, as required.
\end{proof}

To facilitate the presentation of the main result, we provide the following observation (see \cite[Remark 2.1]{HKTT}).

\begin{lem}\label{choose} For any $\alb =(\alpha_1,\ldots, \alpha_n)\in\mathbb R^n$. There exists $\lam=(\lambda_1,\ldots, \lambda_n)\in \mathbb R^n_{+}$ such that $\alb+\lam\in \mathbb \Z^n$ and
$$\sum_{i\in F}\lambda_i < |F|, \text{ for all } F \subseteq \{1,\ldots,n\} \text{ with } F \ne \emptyset.$$
\end{lem}
\begin{proof} Define $\lambda_i = \lceil \alpha_i\rceil-\alpha_i$ for $i=1,\ldots,n$. Then, we can verify that $\lam=(\lambda_1,\ldots, \lambda_n)$ satisfies the lemma.
\end{proof}

We are now in a position to prove the main result of this paper.

\begin{thm} \label{T2} Let $I$ be a square-free monomial ideal in $R$. Then, $\asr I^{(s)} = \asr I^{(s_0)}$ for all $s\geqslant s_0$ where $s_0=\lceil n\bht(I)^{(n+2)/2}\rceil$. In particular, $\asr I^{(s)}$ is constant for $s\gg 0$.
\end{thm}

\begin{proof} To prove the theorem it suffices to show that for any $s\geqslant s_0$ we have $\asr I^{(t)} \subseteq \asr I^{(s)}$ for every $t\geqslant 1$.

We will prove by induction on $n=\dim(R)$. If $n=1$, then $\asr(I^{(s)})$ is constant for every $s\geqslant 1$, and then the theorem holds for this case. Assume that $n\geqslant 2$.

Fix any $t\geqslant 1$. By Lemma \ref{localization} we have
$$\asr(I^{(t)})=\bigcup_{i=1}^n \asr(I[i]^{(t)}) \cup \{J\in \asr(I^{(t)})\mid \supp(J) = [n]\},$$
and
$$\asr(I^{(s)})=\bigcup_{i=1}^n \asr(I[i]^{(s)}) \cup \{J\in \asr(I^{(s)})\mid \supp(J) = [n]\}.$$

For each $i$, we can consider $I[i]$ is a monomial ideal in the polynomial ring $R_i=K[x_j\mid j\ne i]$. Since $\dim(R_i) = n-1$ and $\bht(I[i])\leqslant \bht(I)$, by the induction hypothesis we conclude that $$\asr I[i]^{(t)} \subseteq \asr I[i]^{(s)}.$$

Therefore, it remains to show that $\{J\in \asr(I^{(t)})\mid \supp(J) = [n]\}\subseteq \asr I^{(s)}$. Let $Q\in \asr(I^{(t)})$ with $\supp(Q) = [n]$. Then, we can write $Q=\sqrt{I^{(t)}:\x^{\alb}}$ for some monomial $\x^{\alb}=x_1^{\alpha_1}\ldots x_n^{\alpha_n}\notin I^{(t)}$. 

Suppose that $I$ has the minimal primary decomposition given by $(\ref{square-free-decomp})$. Let $F_i = \supp(\p_i)$ for $i=1,\ldots,r$. We may assume that $Q = \p_1\cap \cdots\cap \p_k$ for some $1\leqslant k\leqslant r$. 

By Lemma \ref{FF}, $\alb$ belongs to the convex polyhedron $\Lambda_t$ in $\R^n$  defined by
\begin{equation}\label{alphaRoot}
\begin{cases}
\sum\limits_{i\in F_j} x_i\leqslant t-1, \text{ for } j=1,\ldots,k,\\
\sum\limits_{i\in F_j} x_i\geqslant t, \text{ for } j=k+1,\ldots,r,\\
x_1\geqslant 0, \ldots,x_n\geqslant 0.
\end{cases}    
\end{equation}
In particular, $\Lambda_t\ne\emptyset$.

Let $\overline{C_1}\subseteq \mathbb R^n$ be the set of solutions of  the  following inequalities system:
\begin{equation}\label{EQ2}
\begin{cases}
\sum\limits_{i\in F_j} x_i  \leqslant  1 & \text{ for } j = 1,\ldots,k,\\
\sum\limits_{i\in F_j} x_i  \geqslant  1 & \text{ for } j =k+1,\ldots,r,\\
x_1\geqslant 0,\ldots,x_n\geqslant 0.
\end{cases}
\end{equation}

Since $\Lambda_t\ne\emptyset$ and $\bigcup_{i=1}^k F_i = \supp(Q) = [n]$, by Lemma \ref{polytope}, $\overline{C_1}$ is a polytope in $\mathbb R^n$ with $\dim \overline{C_1} =n$.

If $\overline{C_1}$ has no supporting hyperplanes of the form $\sum\limits_{i\in F_j} x_i  =1$ for some $j > k$. In this case, we have $k=r$. In particular, $Q = I$, and so $Q\in \asr(I^{(s)})$.

Assume that $\overline{C_1}$ has a supporting hyperplanes of the form $\sum\limits_{i\in F_j} x_i  =  1$ for some $j > k$. We may assume such a supporting hyperplane is $\sum\limits_{i\in F_r} x_i  =  1$;  and let $F$ be
the facet of $\overline{C_1}$ is determined by this hyperplane. Now take $n$ vertices of $\overline{C_1}$ lying in $F$,  namely $\u^1,\ldots,\u^n$, such that they are affinely independent. Let 
$$\u = \frac{1}{n}(\u^1+\cdots+\u^n)$$
be the barycenter of the $(n-1)$-simplex $[\u^1,\ldots,\u^n]$. Then, $\u$ is a relative interior point of $F$, so that it does not belong to any another facet of $\overline{C_1}$.

By Lemma \ref{choose} there is $\lam =(\lambda_1,\ldots, \lambda_n)\in\mathbb R^n_{+}$ such that 
\begin{equation}\label{bound-Lambda}
\sum_{i\in F_j}\lambda_i  < |F_j|, \text{ for all } j=1,\ldots,r,
\end{equation}
and $s\u+\lam\in \mathbb N^n$. Let $\btb =s\u+\lam$. We will show that $Q=\sqrt{I^{(s)}\colon \x^{\btb}}$.

Indeed, suppose $\u = (u_1,\ldots,u_n)$ and $\btb = (\beta_1,\ldots,\beta_n)$. Since $\u\in \overline{C_1}$, for any $j\in \{k+1,\ldots,r\},$ by \eqref{EQ2}, we have
\begin{equation}\label{BETA_ONE}
\sum\limits_{i\in F_j} \beta_i=\sum\limits_{i\in F_j} (su_i+\lambda_i)\geqslant \sum\limits_{i\in F_j} su_i=s\sum\limits_{i\in F_j} u_i\geqslant s,
\end{equation}

For a given $p=1,\ldots, n$, suppose that $\u^p=(u_{p1},\ldots,u_{pn})\in\mathbb R^n$.  As $\u^p$ is a vertex of the polytope $\overline{C_1}$, by \cite[Formula 23 in Page 104]{S}, $\u^p$ is the unique solution of a system of linear equations in the form
\begin{equation*}
\begin{cases}
\sum\limits_{i\in F_j} x_i  =1 & \text{ for } j \in  S_1\subseteq \{1,\ldots,r\},\\
x_i=0 & \text{ for } i \in S_2\subseteq \{1,\ldots,n\},
\end{cases}
\end{equation*}
in which $|S_1|+|S_2|=n.$ By Cramer's rule, we can represent coordinates of $\u^p$ as follows:
\begin{align}\label{equationdet}
u_{pi}=\dfrac{\delta_{pi}}{\delta_p}, \text{ for }i=1,\ldots,n,
\end{align}
where $\delta_p, \delta_{p1},\ldots, \delta_{pn}\in \mathbb N$ and $\delta_p$ is the absolute value of the determinant of this system of linear equations. In particular, $\delta_p\u^p\in \mathbb N^n.$

Let $\delta = \max\{\delta_1,\ldots,\delta_n\}$. Note that the matrix of this system of linear equations has only entries $0$ or $1$, and there are at most ${\rm bight}(I)$ entries $1$ in each row. By Hadamard’s inequality we get
\begin{align}\label{compare}
\delta\leqslant  {\rm bight}(I)^{n/2}.
\end{align}

Fix any $j\in\{1,\ldots,k\}$. Since $\u$  is not lying on the hyperplane $\sum_{i\in F_j}x_i = 1$, so is $\u^m$ for some  $m=1,\ldots,n$. By \eqref{EQ2}, we have
$$\sum\limits_{i\in F_j} u_{mi} < 1.$$
Together with (\ref{equationdet}), it forces
$$\sum\limits_{i\in F_j} u_{mi} \leqslant  1-\frac{1}{\delta}.$$

By combining this inequality with \eqref{EQ2}, we have
\begin{align*}
\sum_{i\in F_j}u_i&= \frac{1}{n}\sum_{p=1}^n \sum_{i\in F_j}u_{pi} = \frac{1}{n}\left( \sum_{i\in F_j}u_{mi}+\sum_{p\ne m} \sum_{i\in F_j}u_{pi}\right)\\ &\leqslant\frac{1}{n} \left(1-\frac{1}{\delta}+\sum_{p\ne m}1\right) =\frac{1}{n} \left(1-\frac{1}{\delta}+n-1\right)=1-\frac{1}{n\delta}. 
\end{align*}
Thus,
\begin{align}\label{compare-u}
\sum_{i\in F_j} u_i \leqslant 1-\frac{1}{n\delta}.
\end{align}

On the other hand, since $s\geqslant s_0\geqslant n{\rm bight} (I)^{(n+2)/2}$, by (\ref{compare}) we get 
\begin{equation}\label{q-bound}
s\geqslant n\delta\bht(I).
\end{equation}

Note that $\bht(I) = \max\{|F_j| \mid j=1,\ldots,r\}$. Combining \eqref{compare-u} and  \eqref{q-bound}, we obtain
\begin{align}\label{compare12}
\sum\limits_{i\in F_j} \beta_i &=\sum\limits_{i\in F_j}(su_i+\lambda_i) =  s\sum\limits_{i\in F_j}u_i+  \sum\limits_{i\in F_j}\lambda_i < s\left(1-\frac{1}{n\delta}\right)+|F_j| \\
&\leqslant s -\bht(I)+|F_j|\leqslant s \notag.
\end{align}
Since the left hand side of \eqref{compare12} is an integer. It follows that 
\begin{align}\label{compare13} 
\sum\limits_{i\in F_j} \beta_i \leqslant s-1.
\end{align}\

Together  $(\ref{BETA_ONE})$ with   \eqref{compare13}, and Lemma \ref{FF}, we conclude that $Q=\sqrt{I^{(s)}:\x^{\btb}}$, so $Q\in \asr(I^{(s)})$. The proof of the theorem is complete.
\end{proof}

We conclude this section with a remark on the asymptotic behavior of the sequence $\{\depth(R/I^{(s)})\}_{s\geqslant 1}$. Recall that the symbolic analytic spread of $I$ is define by
$$\ell_s(I) = \dim \mathcal{R}_s(I)/\mathfrak{m}\mathcal{R}_s(I),$$
where $\mathcal{R}_s=\oplus_{k\geqslant 0} I^{(k)}$ is the symbolic Rees ring of $I$, and $\mathfrak m=(x_1,\ldots,x_n)$ is the maximal homogeneous ideal of $R$. The Theorem \ref{T2} has the following corollary.

\begin{cor}\label{DETP-STABILITY} Let $I$ be a square-free monomial ideal in $R$. Then,
$\depth R/I^{(s)} = \dim R -\ell_s(I)$  for all $s \geqslant n \bht(I)^{(n+2)/2}$.
\end{cor}
\begin{proof} By \cite[Theorem 2.4]{HKTT}
$\depth R/I^{(s)} = \dim R -\ell_s(I)$ for $s\gg 0$. On the other hand, by Hochster's formula \eqref{Hochsterfomular} and Theorem \ref{T2}, we have
$\depth R/I^{(s)}$ is constant for $s\geqslant n\bht(I)^{(n+2)/2}$. Hence,
$\depth R/I^{(s)} = \dim R -\ell_s(I)$ for $s \geqslant n \bht(I)^{(n+2)/2}$, and the corollary follows.
\end{proof}

\section{Powers of cover ideals of balanced hypergraphs}

In this section, we will prove that the set $\asr(I^s)$ is an increasing function in $s$ if $I$ is a cover ideal of a balanced hypergraph.

Let $\H=(\V,\E)$ be a hypergraph on the vertex set $\V=\{1,\ldots,n\}$. A cycle of length $k$ in hypergraph $\H$ is a sequence $(i_1,E_1,i_2,E_2, i_3,\ldots, i_k,E_k,i_1)$ such that: 
\begin{itemize}
\item $E_1,\ldots,E_k$ are distinct edges of $\H$;
\item $i_1,\ldots,i_k$ are distinct vertices of $\H$ such that $\{i_t,i_{t+1}\}\subseteq E_t$ with $t=1,\ldots,k-1$;
\item $\{i_k,i_1\}\subseteq E_k.$
\end{itemize}

A hypergraph $\H$ is said to be  balanced if every odd cycle has an edge that contains at least three vertices of the cycle.

This notion also can be defined by using incidence matrices. The incidence matrix $A(\H) = (a_{ij})_{m\times n}$, with rows representing the edges $E_1,E_2, \ldots,E_m$, and  columns representing the vertices $1,2,\ldots, n$ where $a_{ij} = 0$ if $j\notin E_i,$ and  $a_{ij} = 1$ if $j\in E_i$.

A matrix is called a balanced matrix if it 
 has no square submatrix of the form
 \[B_k= \begin{pmatrix} 1&1&0& \cdots &0&0 & 0 \\
0&1&1&\cdots &0& 0 &0  \\
0&0&1&\cdots &0&0&0 \\
\vdots &\vdots &\vdots & \cdots & \vdots & \vdots& \vdots\\

0 &0 &0 & \cdots & 0&1 & 1\\
1 &0 &0 & \cdots &0& 0 & 1
\end{pmatrix},\]
 where $k\geq 3$ is odd.

With this definition, $\H$ is balanced if and only if $A(\H)$ is balanced.

\begin{lem}\label{balanced} \cite[Corollary $1.6$]{HHTZ} If $\H$ is balanced, then $I(\H)^s = I(\H)^{(s)}$ for all $s\geqslant 1$.
\end{lem}

We now prove the main result of this section.

\begin{thm} \label{T3} Let $I=J(\H)$ be a cover ideal of balanced hypergraph $\H$. Then, $\asr (I^s)\subseteq \asr (I^{s+1})$ for every $s\geqslant 1$.
\end{thm}

\begin{proof} By Lemma \ref{balanced}, it suffices to show that $\asr(I^{(s)}) \subseteq \asr(I^{(s+1)})$ for every $s\geqslant 1$.

Suppose that $\E = \{E_1,\ldots,E_m\}$ be the set of edges of $\H$. By (\ref{square-free-dmp}), we have 
$$I = \bigcap_{j=1}^m \p_j$$
where $\p_j=(x_i |\ i\in E_j)$, for $j=1,\ldots,m$.

We now prove $\asr(I^{(s)})\subseteq \asr(I^{(s+1)})$ by induction on $n$. If $n=1$, then we can verify that $\asr(I^{(s)})=\{(x_1)\}$ for all $s\geqslant 1$. Hence, the theorem holds for this case.

Assume that $n\geqslant 2$. For any $s\geqslant 1$, and $Q\in \asr(I^{(s)})$, we need to show that $Q\in\asr(I^{(s+1)})$.

By Lemma \ref{localization} we have
$$\asr(I^{(s)})=\bigcup_{i=1}^n \asr(I[i]^{(s)}) \cup \{J\in \asr(I^{(s)})\mid \supp(J) = [n]\},$$
and
$$\asr(I^{(s+1)})=\bigcup_{i=1}^n \asr(I[i]^{(s+1)}) \cup \{J\in \asr(I^{(s+1)})\mid \supp(J) = [n]\}.$$

For each $i$, we can consider $I[i]$ is a monomial ideal in the polynomial ring $R_i=K[x_j\mid j\ne i]$. Let $\H_i$ be the hypergraph with the edge $\{E\mid E\in \E \text{ and } i\notin E\}$. Then, $\H_i$ is also a balanced hypergraph. Observe that $I[i]$ is the cover ideal of $\H_i$. On the other hand, since $\dim(R_i) = n-1$ and $\bht(I[i])\leqslant \bht(I)$, by the induction hypothesis we conclude that $$\asr (I[i]^{(s)}) \subseteq \asr (I[i]^{(s+1)}).$$

Therefore, it remains to show that $\{J\in \asr(I^{(s)})\mid \supp(J) = [n]\}\subseteq \asr I^{(s+1)}$. Let $Q\in \asr(I^{(s)})$ with $\supp(Q) = [n]$. Then, we can write $Q=\sqrt{I^{(s)}:\x^{\alb}}$ for some monomial $\x^{\alb}=x_1^{\alpha_1}\ldots x_n^{\alpha_n}\notin I^{(s)}$. 

We may assume that $Q = \p_1\cap \cdots\cap \p_k$ for some $1\leqslant k\leqslant m$. Let $F_i = \supp(\p_i)$ for $i=1,\ldots,r$.

By Lemma \ref{FF},  $\alb$ belongs to the convex polyhedron $\Lambda_s$ in $\R^n$ defined by the linear system:
\begin{equation} \label{solPs}
\begin{cases}
\sum\limits_{i\in F_j} x_i \leqslant s-1 & \text{ for } j = 1,\ldots,k,\\
\sum\limits_{i\in F_j} x_i \geqslant  s & \text{ for } j = k+1,\ldots,m,\\
x_1\geqslant 0, \ldots,x_n\geqslant 0.
\end{cases}
\end{equation}

Next, let $\overline{C_1}$ be the polyhedron in $\R^n$ defined by
\begin{equation}\label{EQ-basics_1}
\begin{cases}
\sum\limits_{i\in F_j} x_i  \leqslant  1 & \text{ for } j =1,\ldots,k,\\
\sum\limits_{i\in F_j} x_i  \geqslant  1 & \text{ for } j =k+1,\ldots,m,\\
x_1\geqslant 0,\ldots,x_n\geqslant 0.
\end{cases}
\end{equation}

Since $\Lambda_s\ne \emptyset$ and $\supp(Q)=\bigcup_{i=1}^k F_i = [n]$, by Lemma \ref{polytope} we have $\overline{C_1}$ is a polytope with $\dim(\overline{C_1})=n$. Let $\gmb\in \mathbb R^n$ be a vertex of $\overline{C_1}$. Then, $\gmb\in\N^n$ by \cite[Lemma 2.1]{HG}.

Let $\btb =\alb+\gmb\in \N^n$. Since $\alb$ is a solution of the linear system $(\ref{solPs})$ and $\gmb$ is a solution of the linear system $(\ref{EQ-basics_1})$, we can verify that
\begin{equation}\label{mut+1}
\begin{cases}
\sum\limits_{i\in F_j} \beta_i \leqslant s & \text{ for } j = 1,\ldots,k,\\
\sum\limits_{i\in F_j} \beta _i  \geqslant  s+1 & \text{ for } j = k+1,\ldots,m.
\end{cases}
\end{equation}

Together with Lemma \ref{FF},  it yields
$$\sqrt{I^{(s+1)}:\x^{\btb}}=\pp_1\cap \ldots \cap \pp_k=Q,$$
so $Q\in \asr (I^{(s+1)})$, and the theorem follows.
\end{proof}

\begin{cor} Let $G$ be a bipartite graph. Then, $\asr(J(G)^s)$ is increasing in $s$.
\end{cor}
\begin{proof} Since $G$ is bipartite, it has no odd cycle, so it is balanced. Thus, the corollary follows from Theorem \ref{T3}.
\end{proof}

According to this corollary, it is natural to ask whether $\asr(J(G)^s)$ is increasing in $s$ for every graph $G$. It is worth mentioning that Francisco, H{$\rm \grave{a}$} and Van Tuyl conjectured that $\ass(J(G)^s)$ is increasing in $s$ (see \cite{FHV}), but a counter-example is given in \cite{KSS}. By using this counter-example, we can prove the following proposition.

\begin{prop} There is a graph $G$ such that $\asr(J(G)^s)$ is not increasing in $s$.
\end{prop}
\begin{proof} By \cite[Theorem 11]{KSS}, there is a graph $G$ such that $\ass(J(G)^s)$ is not increasing in $s$. Together with Lemma \ref{inc-lemma}, it implies that $\asr(J(G)^s)$ is also not increasing in $s$.
\end{proof}

\subsection*{Acknowledgment}  Truong Thi Hien was funded by the Postdoctoral Scholarship Programme of Vingroup Innovation Foundation (VINIF), code [VINIF.2024.STS.34].

\end{document}